\documentclass[11pt]{amsart}
\usepackage[normalem]{ulem}
\usepackage[usenames]{color}
\usepackage{graphicx,amscd,psfrag,amssymb}
\usepackage[colorlinks=true,citecolor=red,linkcolor=blue]{hyperref}

\DeclareMathOperator{\area}{area}

\DeclareMathOperator{\arc}{arc}
\DeclareMathOperator{\sector}{sector}
\DeclareMathOperator{\quadrilateral}{quadrilateral}
\DeclareMathOperator{\perimeter}{Perimeter}

\theoremstyle{plain}

\theoremstyle{definition}

\theoremstyle{remark}

\newcommand{\ncite}[2][]{\nocite{#2}}

\begin{document}

\title{Circular reasoning: who first proved that $C/d$ is a constant?}
\author{David Richeson}   \address{Dickinson College\\ Carlisle, PA 17013} \email{richesod@dickinson.edu} 

\maketitle

\begin{abstract}
We answer the question: who first proved that $C/d$ is a constant? We argue that Archimedes proved that the ratio of the circumference of a circle to its diameter is a constant independent of the circle and that the circumference constant equals the area constant ($C/d=A/r^{2}$). He stated neither result explicitly, but both are implied by his work. His proof required the addition of two axioms beyond those in Euclid's \emph{Elements}; this was the first step toward a rigorous theory of arc length. We also discuss how Archimedes's work coexisted with the 2000-year belief---championed by scholars from Aristotle to Descartes---that it is impossible to find the ratio of a curved line to a straight line.
\end{abstract}

\ncite{Seidenberg:1972}
\ncite{Boyer:1964}
\ncite{Dijksterhuis:1987}
\ncite{Traub:1984}
\ncite{Richman:1993}
\ncite{Gillman:1991}

For a long time I was too embarrassed to ask the question in the title. The expression $C=\pi d$, which gives the relationship between the circumference and the diameter of a circle, is one of the few formulas known to almost all children and adults, regardless of how long they have been out of school. But who first proved it?  I was sure that merely asking the question would tarnish my reputation as a mathematician and as a student of the history of mathematics. Surely the answer is well-known and easy to find. But I could not locate it---at least none of the standard reference books discussed the matter. 

When I finally summoned the courage to ask others I received two types of replies: either ``It is obvious; all circles are similar'' or ``I don't know.'' The first of these is intriguing and has some merit; it may be why the invariance of $C/d$ was discovered in so many different cultures. But it does not easily lead to a rigorous proof.

So I began to dig deeper, and the more I looked, the more confusing were the facts. In short:

\begin{itemize}
\item
The existence of the circumference constant ($C/d$) was known to many ancient civilizations, and thus one must expect it was known to the Greeks.
\item
Aristotle (384--322 BCE) asserted that it was impossible to compare the lengths of curves and line segments, and this belief was widely held until the seventeenth century. Thus one cannot discuss the ratio of the circumference to the diameter.
\item
In \emph{Elements} Euclid (fl.\ c.\ 300 BCE) proved, essentially, the existence of the area constant for circles ($A/r^{2}$). But he did not mention the invariance of $C/d$ or anything equivalent to it.
\item
Archimedes (287--212 BCE) proved that $A=\frac{1}{2}Cr$ and he showed that $223/71<C/d<22/7$. But he never explicitly stated that $C/d$ is a constant. 
\end{itemize}

After a deeper immersion in the literature I discovered the answer to my question: Archimedes. In this article we describe how Archimedes extended Euclid's axioms so that he could treat arc length with the rigor demanded by the Greeks.  His work on arc length has three important, \emph{unstated} corollaries: $C/d$ is a constant, the circumference constant equals the area constant, and the problem of rectifying the circle is equivalent to the problem of squaring the circle. The first of these is the theorem that we are seeking; Archimedes may have stated the result explicitly, but no record of that remains.

We also examine the long-held belief---which sounds ridiculously counterintuitive to us---that curves are not like strings that can be straightened into lines. That is, lines and curves are fundamentally different and cannot be compared; we cannot compute the ratio of a curve to a segment. 
We recount the two thousand year discussion of the difference between straight and curved lines and the definitions of arc length. Archimedes's theorem is an important part of this conversation. Despite the fact that he set a line segment (the side of a triangle) equal to the circumference of a circle, the resulting theorem reinforced Aristotle's assertion that lines and curves cannot be compared.  It shows that rectifying the circle is equivalent to squaring the circle, and because the problem of squaring the circle was believed impossible, so must be the problem of rectifying the circle.

We must point out that throughout this article we will use modern terminology, notation, and mathematical concepts unless it is important to understand the original approach. The phrase ``the value $A/r^{2}$ is independent of the circle and it is equal to the circumference constant $C/d=\pi$''  would, for many reasons, be completely foreign to the Greeks. For example, they would express these ideas through ratios and proportions, and the terms in these ratios, $A$, $r^{2}$, $C$, and $d$, would  be geometric objects, not numbers. And of course the use of $\pi$ to represent this constant was many years in the future. (The first use is attributed to William Jones (1675--1749) in 1706 \cite{Jones:1706}.)

\section{Is there something to prove?}

One common response to the question in the title is: ``It is obvious. All circles are similar.'' For example, Katz writes that ``It may be obvious that the circumference is proportional to the diameter\ldots'' \cite[p.~20]{Katz:1998} 
Indeed, if we scale a figure by a factor of $k$, then all lengths increase by this same factor. Thus the ratios of corresponding diameters and circumferences remain constant. 
The ``it's obvious'' argument probably explains why it was rediscovered by so many cultures. 

But this response is mathematically unsatisfying; it is not a simple processes to turn it into a rigorous proof, and Euclid did not attempt to do so. He defined similarity for polygons (Def.~VI.1: ``Similar rectilinear figures are such as have their angles severally equal and the sides about the equal angles proportional.'' \cite[p.~188]{Heath:1908a}) but not for circles or other curves.  It is not clear how he would do so. He could have said that two circles are similar provided $C_{1}:d_{1}::C_{2}:d_{2}$, but of course that is what we want to prove.





Today we define similarity in terms of functions; a function $f$ from the Euclidean plane to itself is a \emph{similarity transformation} if there is a positive number $k$ such that for any points $x$ and $y$, $d(f(x),f(y))=k\cdot d(x,y)$. Two geometric shapes are \emph{similar} if one is the image of the other under such a mapping. If $AB$ and $CD$ are two segments in a geometric figure and $A'B'$ and $C'D'$ are the corresponding segments in a similar figure, then, assuming that the length of a line segment is the distance between its endpoints, \[\frac{A'B'}{C'D'}=\frac{d(A',B')}{d(C',D')}=\frac{d(f(A),f(B))}{d(f(C),f(D))}=\frac{k\cdot d(A,B)}{k\cdot d(C,D)}=\frac{d(A,B)}{d(C,D)}=\frac{AB}{CD}.\]  Despite the fact that  it is easy to show that two circles in the plane are similar, we still cannot say that $C/d$ is constant. To do so we need a definition for the length of the circumference; but it is not clear how to extend the definition for the length of a segment to the length of a curve. This brings us to the key question: what \emph{is} the length of a curve? How do we define arc length? 

The history of arc length is long and fascinating. For a thorough discussion see Gilbert Traub's Ph.D. dissertation \cite{Traub:1984}. In this article we are most interested in the length of the circumference of a circle, which is a much easier question to answer than the same question for general curves. We say a little about arc lengths of other curves at the end of the article.

\section{The early history of $\pi$}

Books and articles on the history of mathematics are quick to point out that ancient civilizations knew about the circle constant $\pi$---or circle constants, we should say, for $\pi$ plays a dual role as the area constant and the circumference constant. A great deal has been written about whether the Egyptians, the Babylonians, the Chinese, the Indians, and the writers of the Bible knew about the constants, about what value they used for them, and about whether they knew that the two constants were the same. Thus one cannot help but imagine the Greeks knew of this constant and were looking for rigorous proofs of its existence. (For a history of $\pi$ see, for example, \cite{Beckmann:1971,Berggren:2004,Schepler:1950a,Schepler:1950b,Schepler:1950c,Castellanos:1988a,Castellanos:1988}.)

Proposition XII.2 of Euclid's \emph{Elements} implies the existence of the area constant $\pi$. It states ``Circles are to one another as the squares on their diameters.'' \cite[p.~371]{Heath:1908b} 
Symbolically, if we have circles with areas $A_{1}$ and $A_{2}$ and diameters $d_{1}$ and $d_{2}$, respectively, then \[\frac{A_{1}}{A_{2}}=\frac{d_{1}^{2}}{d_{2}^{2}}.\] In other words, $A/d^{2}$ is the same value for every circle. Thus so is $A/r^{2}$. This theorem may have been proved first by Hippocrates of Chios  (c.\ 470--c.\ 410 BCE), but Euclid's proof---a classic application of the method of exhaustion---is due to Eudoxus of Cnidus (c.\ 400--c.\ 347 BCE). 

One expects to see a proposition in Euclid's \emph{Elements} such as:  ``The circumferences of circles are to one another as their diameters.'' But on this topic Euclid is silent. In fact, Euclid never compared a line segment to the arc of a circle (or any other curve). 
He only compared two magnitudes (lengths, areas, angles, etc.) if they are of the same type. In Definition V.4, he wrote that  ``Magnitudes are said to have a ratio to one another which can, when multiplied, exceed one another.'' \cite[p.~114]{Heath:1908a} Thus it must be possible to put together several copies of one magnitude to obtain a magnitude larger than the other, and vice versa. This definition was flexible enough for him to find ratios of incommensurable magnitudes (such as the ratio of the diagonal of a square to its side).  But it did not allow him to find the ratio of an angle to a line segment, an area to a volume, or---and this is the key fact for us---a circular arc to a line segment.\footnote{He also did not compute the ratio of two arcs of different radii, but comparing arcs of circles with the same radius is fine. See Proposition VI.33 and the discussion that follows it in \cite[pp.~273--276]{Heath:1908a}.}

The belief that straight lines and curves are fundamentally different---and are not comparable---is attributed to Aristotle (384--322 BCE). In his \emph{Physics} he wrote about whether it is possible to compare the motion of a body that travels on a circular path with a body that moves in a straight line. He concluded that it is not. 
\begin{quote}
The fact remains that if the motions are comparable, there will be a straight line equal to a circle. But the lines are not comparable; so neither are the motions. 
\cite[p.~426]{Ross:1936}, \cite[pp.~140--142]{Heath:1980}.
\end{quote}
To put this more concretely, the Greeks believed that it was not possible to construct, using only a compass and straightedge, a line segment equal in length to a given curve. That is, they believed that it is impossible to rectify a curve using Euclidean tools. 

\section{Euclid's common notion 5}

The reason that the theorem on circular areas was proved before the theorem on circumferences is that arc length is inherently more complicated than area (just ask any first-year calculus student!). Gottfried Leibniz (1646--1716) wrote that ``areas are more easily dealt with than curves, because they can be cut up and resolved in more ways.'' \cite[p.~61]{Traub:1984} In particular, Euclid could apply his fifth common notion---``The whole is greater than the part'' \cite[p.~155]{Heath:1908}---to areas, but not to curves.\footnote{Proclus (411--485) claimed that this common notion was in Euclid's original version of \emph{Elements}, whereas Heath suspected that it was added later \cite[p.~232]{Heath:1908}. Either way, whether it was explicit or implicit in \emph{Elements}, the argument here still holds.} For example, suppose a circle is inscribed in an equilateral triangle. The circle is a subset of the triangle. By Euclid's common notion, the area of the circle is less than the area of the triangle. However, while it is intuitively clear that the circumference of the circle is shorter than the perimeter of the triangle, Euclid's common notion does not guarantee this.

\begin{figure}[ht]
\includegraphics{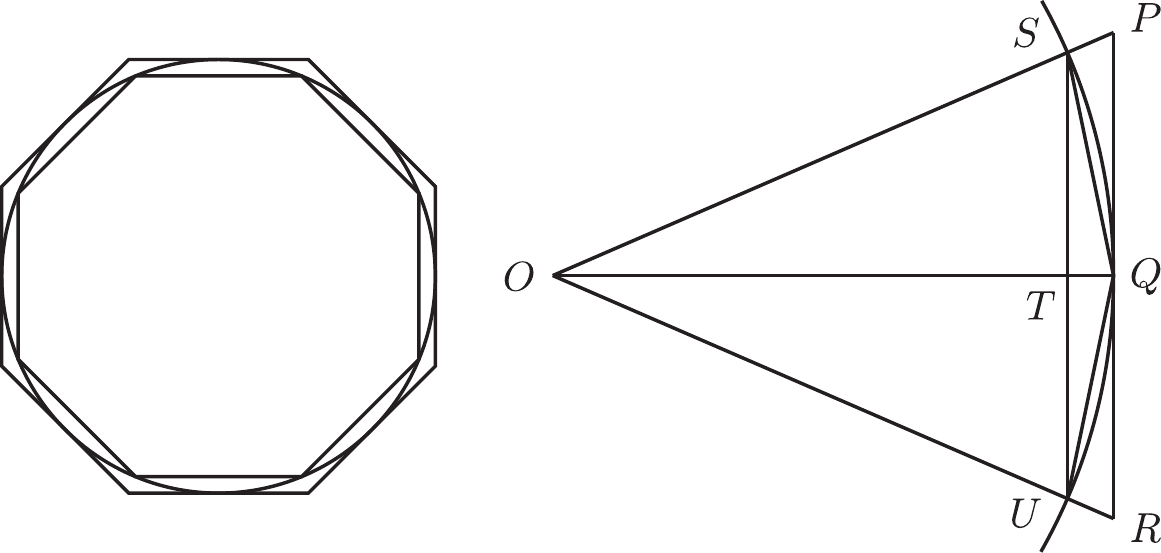}
\caption{}
\label{fig:sector}
\end{figure}

As an exercise, let us try to prove that the perimeter of an inscribed regular $n$-gon is less than the circumference of the circle, which is less than the perimeter of a circumscribed regular $n$-gon (see Figure \ref{fig:sector}). Suppose the circle has radius $r$, circumference $C$, and area $A$. For simplicity we may restrict our attention to a sector of the circle in which $SU$ and $PR$ are sides of the inscribed and circumscribed polygons, respectively.  We would like to prove that $SU<\arc(SQU)<PR$.

We know that
\[\quadrilateral OSQU=\triangle OSQ\cup\triangle OQU\subset\sector OSQU\subset \triangle OPR.\]
By Euclid's fifth common notion
\[\area(\triangle OSQ)+\area(\triangle OQU)<\area(\sector OSQU)<\area(\triangle OPR),\]
and hence
\[\frac{1}{2}OQ\cdot ST+\frac{1}{2}OQ\cdot TU<\area(\sector OSQU)<\frac{1}{2}OQ\cdot PR.\]
We know that $OQ=r$ and $ST+TU=SU$, so
\[\frac{1}{2}SU\cdot r<\area(\sector OSQU)<\frac{1}{2}PR\cdot r,\]
and 
\[SU<\frac{2\area(\sector OSQU)}{r}<PR.\]
Because 
\[\frac{\area(\sector OSQU)}{A}=\frac{\arc(SQU)}{C},\]
it follows that
\[SU<\frac{2A}{Cr}\arc(SQU)<PR.\]
To conclude that 
\[SU<\arc(SQU)<PR\]
we need the fact that $2A/Cr=1$. But this is precisely Archimedes's theorem, and, as we shall see, his proof used the fact that \[\perimeter(\text{inscribed polygon})<C<\perimeter(\text{circumscribed polygon}).\] The argument is circular (pun intended)!\ncite{Seidenberg:1972}\ncite{Richman:1993}\footnote{This conversation is intimately related to that of how to prove that $\sin\theta<\theta<\tan\theta$, a key ingredient in the proof that $\displaystyle\lim_{\theta\to 0}\frac{\sin\theta}{\theta}=1$. In our diagram, if $\theta=\angle SOQ$, then $SU=2r\sin\theta$, $PR=2r\tan\theta$, and---assuming that angle measures arc length---$\arc (SQU)=2r\theta$. (See \cite{Richman:1993}, \cite[pp.~216--220]{Comenetz:2002}, and \cite{Gillman:1991} for more on this topic.)}

One must suspect that even the early Greeks were aware of the relationship between the circumference and diameter. And one must imagine that Euclid would have loved nothing more than to include the ``circumference theorem'' in \emph{Elements}. But he was unable to give a rigorous proof of it using his postulates and common notions.  Moreover, the mathematical wisdom of the day said that proving it was impossible. Indeed, it took the genius of Archimedes to realize that to prove this theorem we need to add two axioms to those given by Euclid.\ncite[p.~36]{Traub:1984}

\section{Measurement of a Circle}

To assemble Archimedes's proof that $C/d$ is a constant we need two of his works: \emph{Measurment of a Circle} \cite[pp.~91--98]{Archimedes:2002} and \emph{On the Sphere and Cylinder I} \cite[pp.~1--55]{Archimedes:2002}. The first contains the proof and the second contains the new axioms that are required for the proof. The order in which they were written is unknown. Historians once believed that \emph{On the Sphere and Cylinder} came first because that makes the most sense logically---axioms first, theorem second (see \cite{Archimedes:2002}, for instance). However, more recently scholars have argued that they were written in the other order. Knorr speculated that when writing the \emph{Measurement of a Circle} Archimedes accepted the axioms as obvious truths.  ``The formulation of these principles as explicit axioms in \emph{On the Sphere and Cylinder} would thus result from Archimedes' own later reflections on the formal requirements of such demonstrations.'' \cite[pp.~153--155]{Knorr:1993}
 

The surviving copy of \emph{Measurment of a Circle} has only three propositions and is not a faithful copy of Archimedes's writing (one of the three results is clearly incorrect as stated). Dijksterhuis refers to it as ``scrappy and rather careless''  and points out that ``it is quite possible that the fragment we possess formed part of a longer work, which is quoted by Pappus under the title \emph{On the Circumference of the Circle}.'' \cite[p.~222]{Dijksterhuis:1987}

\begin{figure}[ht]
\includegraphics{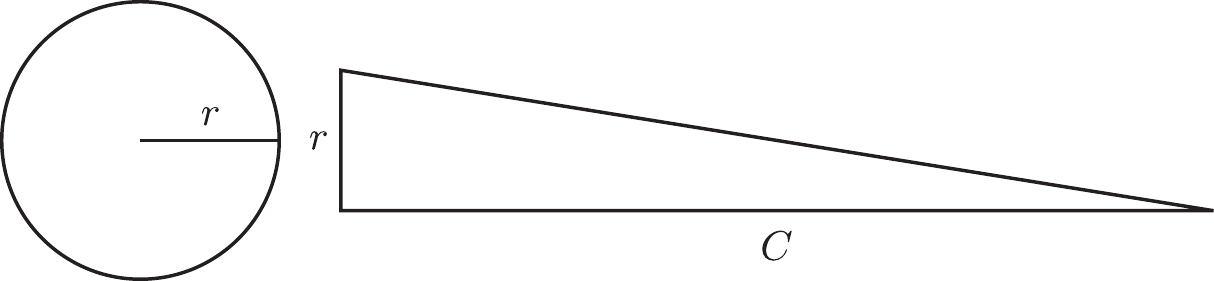}
\caption{}
\label{fig:areacircumference}
\end{figure}

In his first proposition Archimedes broke the Aristotelian rule about not comparing curves and lines. He proved (see Figure \ref{fig:areacircumference}):
\begin{quote}
The area of any circle is equal to a right-angled triangle in which one of the sides about the right angle is equal to the radius, and the other to the circumference, of the circle.\footnote{Knorr argued that the original theorem proved by Archimedes did not create a triangle, but a rectangle. He claimed that the original theorem read something like ``The product of the perimeter of the circle and the (line) from its center is double of the area of the circle.'' He used ``the product'' as an abbreviation for ``the rectangle bounded by.'' \cite{Knorr:1986}}  \cite[p.~91]{Archimedes:2002}
\end{quote}
That is, \[A=\frac{1}{2}Cr.\] Clearly, he was making a line segment (a side of a triangle) equal to the circumference of a circle.

This deliberate decision was not lost on future mathematicians. Theon of Alexandria (c.\ 335--c.\  405), in a commentary on Ptolemy's (85--165) \emph{Collection} Book VI, wrote ``It has been proved by Archimedes that when the periphery of the circle has been unrolled into a straight line\ldots'' \cite{Knorr:1986}. Some emended versions of \emph{Measurement of a Circle} published in the middle ages also addressed this point. An early version, 
which may have been written in the thirteenth century, added three postulates, one of which was that ``There is some curved line equal to any straight line and some straight line equal to any curved line.'' \cite[p.~69]{Clagett:1964} Another version (which can only be placed between the middle of the thirteenth and fifteenth centuries) said
\begin{quote}
The second of the postulates is that a curved line is equal to a straight line. We postulate this although it is a principle known \emph{per se} and recognized by anybody with a sound head. For if a hair or silk thread is bent around circumference-wise in a plane surface and then afterwards is extended in a straight line in the same plane, who will doubt---unless he is hare-brained---that the hair or thread is the same whether it is bent circumference-wise or extended in a straight line and is just as long the other time as the other. \cite[p.~171]{Clagett:1964}
\end{quote}

Some historians argue that Archimedes's theorem relating the circumference and the area was proved one to two hundred years before Archimedes because it was used to exhibit the circle squaring properties of the curve now known as the quadratrix. But more recent scholarship shows that this use of the quadratrix was not due to Hippias of Elis (born c.\ 460 BCE) who discovered the curve and used it to trisect angles, or to Dinostratus (fl.\ c.\ 350 BCE) as others contend, but to Nicomedes (fl.\ c.\ 250 BCE), who was able to benefit from the work of Archimedes and Eudoxus (see \cite[pp.~80--82, 233]{Knorr:1993} or \cite{Smeur:1970} for details). We say more about this curve shortly.

With Archimedes's result and Euclid's XII.2 we are but a few algebraic steps away from proving that $C/d$ is a constant and that it is the same as the area constant: \[\frac{C}{d}=\left(\frac{2A}{r}\right)\frac{1}{d}=\frac{A}{r^{2}}=\pi.\] 

Although Archimedes never stated this invariance of $C/d$ explicitly, it is implied by the third proposition, 
which gives his famous bounds: 
 \[223/71<C/d<22/7.\]
Many people highlight this pair of inequalities as the key result of this work. One popular textbook calls it ``the most important proposition in \emph{Measurement of the Circle.}'' \cite[p.~202]{Burton:2007} Indeed it is amazing that Archimedes overcame the clumsy Greek numerical system to obtain these bounds and that he devised an algorithm that was used for centuries by digit-hunters searching for increasingly more accurate bounds for $\pi$ (approximating the circle with inscribed and circumscribed polygons of many sides). But we argue that the first theorem is the more important one, for it shows that $C/d$ is constant and that this circumference constant equals the area constant. 

\begin{figure}[ht]
\includegraphics{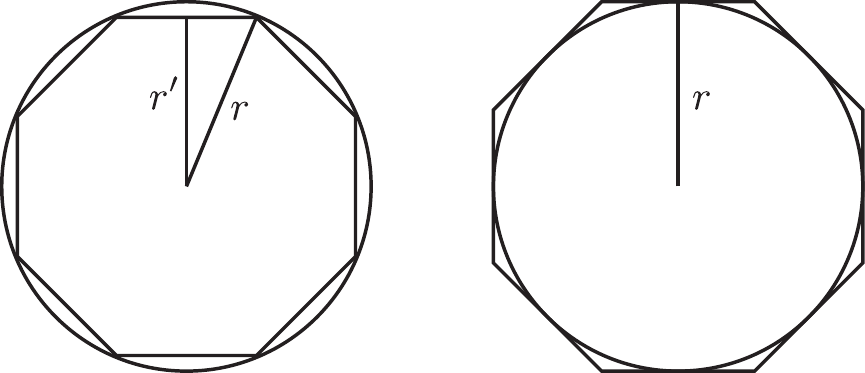}
\caption{}
\label{fig:archimedesproof}
\end{figure}

We now give a sketch of Archimedes's proof that $A=\frac{1}{2}Cr$. Suppose we begin with a circle with radius $r$, area $A$, and circumference $C$. Let $T$ be a right triangle with side lengths $r$ and $C$. For the sake of contradiction, suppose $\area(T)\ne A$.  There are two cases, either $\area(T)< A$ or $\area(T)> A$. Suppose $\area(T)< A$; that is $A-\area(T)>0$. Then Archimedes proved that there is an inscribed regular polygon $P_{\text{in}}$ such that $A-\area(P_{\text{in}})<A-\area(T)$. So $\area(P_{\text{in}})>\area(T)$. Let $r'$ be the length of the apothem of $P_{\text{in}}$. Then $r'<r$ (see Figure \ref{fig:archimedesproof}). It follows that \[\area(P_{\text{in}})=\frac{1}{2}(r')(\perimeter(P_{\text{in}}))<\frac{1}{2}rC=\area(T).\] This is a contradiction. 

Now suppose $\area(T)> A$. Again, Archimedes proved that there is a circumscribed regular polygon $P_{\text{circ}}$ such that $\area(P_{\text{circ}})-A<\area(T)-A$. So $\area(P_{\text{circ}})<\area(T)$. In this case the apothem of $P_{\text{circ}}$ is the radius of the circle. So we have \[\area(P_{\text{circ}})=\frac{1}{2}(r)(\perimeter(P_{\text{in}}))>\frac{1}{2}rC=\area(T),\] a contradiction.

\begin{figure}[ht]
\includegraphics{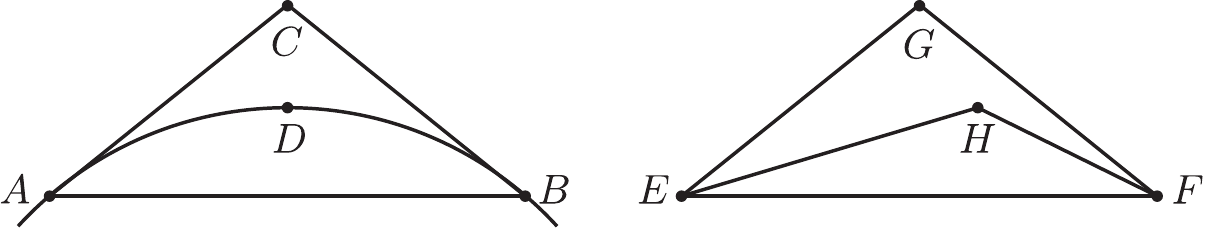}
\caption{}
\label{fig:archimedes2}
\end{figure}

Archimedes used the inequalities \[\perimeter(P_{\text{in}})<C<\perimeter(P_{\text{circ}}),\] without proof. These inequalities would follow if we knew that \[AB<ADB<AC+BC,\] where $AB$ is a chord in a circle, $ADB$ is the included arc, and $AC$ and $BC$ are two segments tangent to the circle (as in Figure \ref{fig:archimedes2}). If we believe Knorr's chronology, it is likely Archimedes assumed these inequalities when writing \emph{Measurement of a Circle}, only later returning to it with a more rigorous mindset.

One may argue that these inequalities are obvious, especially the first, because a line is the shortest distance between two points. Yet Euclid did not assume this as an axiom. In fact he proved both inequalities for segments, namely, that in Figure \ref{fig:archimedes2}, $EF<EH+FH<EG+FG$.  Proposition I.20 is the triangle inequality.
\begin{quote}
In any triangle two sides taken together in any manner are greater than the remaining one. \cite[p.~286]{Heath:1908}
\end{quote}
From this we conclude that $EF<EH+FH$. According to Proclus, the Epicureans scoffed at this theorem, saying that it was so intuitive that even an ass knew it was true: if you put food at one vertex of a triangle and an ass at another it would certainly walk along the edge between them and not along the other two \cite[p.~287]{Heath:1908}. Euclid's next proposition implies that $EG+FG>EH+FH$. 
\begin{quote}
If on one of the sides of a triangle, from its extremities, there be constructed two straight lines meeting within the triangle, the straight lines so constructed will be less than the remaining two sides of the triangle, but will contain a greater angle. \cite[p.~289]{Heath:1908}
\end{quote}

\section{On the Sphere and Cylinder I}
\ncite[pp.~141--154]{Dijksterhuis:1987}

\begin{figure}[ht]
\includegraphics{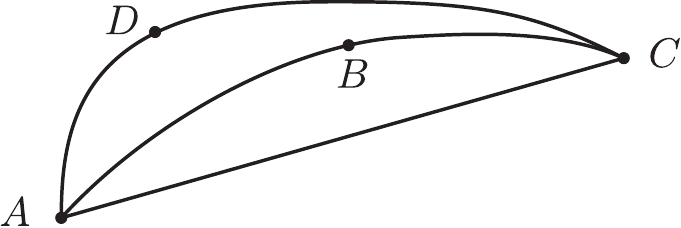}
\caption{}
\label{fig:archimedesaxioms}
\end{figure}

Archimedes began \emph{On the Sphere and Cylinder I} by formalizing the notions that he took for granted in \emph{Measurement of a Circle.} He stated that a curve is \emph{concave in the same direction} 
\begin{quote}
if any two points on it are taken, either all the straight lines connecting the points fall on the same side of the line, or some fall on one and the same side while others fall on the line itself, but none on the other side. \cite[pp.~2]{Archimedes:2002}
\end{quote} 
For example, the curve $ABC$ in Figure \ref{fig:archimedesaxioms} is concave in the same direction (as are $ADC$ and $AC$). Then he stated the following axioms \cite[pp.~3--4]{Archimedes:2002}.
\begin{enumerate}
\item
Of all lines which have the same extremities the straight line is the least.
\item
Of other lines in a plane and having the same extremities, [any two] such are unequal whenever both are concave in the same direction and one of them is either wholly included between the other and the straight line which has the same extremities with it, or is partly included by, and is partly common with, the other; and that [line] which is included is the lesser [of the two].
\end{enumerate}

Returning to Figure \ref{fig:archimedesaxioms}, Archimedes's axiom (1) implies that $AC$ is shorter than $ABC$ and $ADC$, and (2) implies that $ABC$ is shorter than $ADC$. Clearly these axioms are the analogues of Euclid's propositions for piecewise-linear curves.\ncite[p.~34]{Traub:1984} With these axioms in place, Archimedes was ready to make rigorous statements about lengths of curves.  (Archimedes should have included a third axiom: finite additivity. He used the fact that if a curve is broken into parts, then the total length equals the sum of the individual lengths. Euclid's \emph{Elements} should also have contained this axiom. To make up for this Christopher Clavius (1538--1612) added the axiom ``the whole is equal to the sum of its parts'' to his 1574 version \emph{Elements} \cite[p.~40]{Traub:1984}, \cite[p.~323]{Heath:1908}.)

Following the statements of the axioms, Archimedes immediately applied the first one. 
\begin{quote}
If a polygon be inscribed in a circle, it is plain that the perimeter of the inscribed polygon is less than the circumference of the circle; for each of the sides of the polygon is less than that part of the circumference of the circle which is cut off by it. \cite[pp.~4]{Archimedes:2002} 
\end{quote}
Then he stated his first proposition and used the second axiom to prove it.  
\begin{quote}
If a polygon be circumscribed about a circle, the perimeter of the circumscribed polygon is greater than the perimeter of the circle. \cite[pp.~5]{Archimedes:2002} 
\end{quote}

Thus, putting together the proof from \emph{Measurement of a Circle} and the axioms from \emph{On the Sphere and Cylinder} (regardless of the order of their publication), we obtain the theorem that $C/d$ is a constant.



\section{The comparability of lines and curves}

\begin{figure}[ht]
\includegraphics{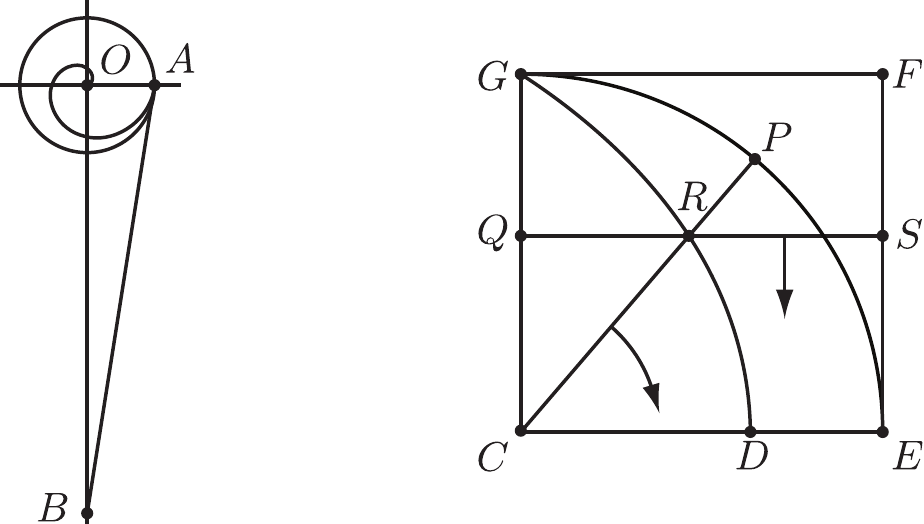}
\caption{}
\label{fig:quadratrixspiral}
\end{figure}

Archimedes's theorem relating the area of a circle to the area of a triangle was not the only Greek instance in which a curve was shown to be or set to be equal to a straight line. In \emph{On Spirals} \cite[pp.~151--188]{Archimedes:2002} Archimedes gave a method of rectifying a circle that used the arithmetic (now Archimedean) spiral (see Figure \ref{fig:quadratrixspiral}).  Proposition 18(I) states 
\begin{quote}
If $OA$ be the initial line, $A$ the end of the first turn of the spiral, and if the tangent to the spiral at $A$ be drawn, the straight line $OB$ drawn from $O$ perpendicular to $OA$ will meet the said tangent in some point $B$, and $OB$ will be equal to the circumference of the `first circle.'\,'' \cite[p.~171]{Archimedes:2002}
\end{quote}

Also, Nicomedes used the quadratrix to rectify the circle. The quadratrix is drawn as follows (see Figure \ref{fig:quadratrixspiral}). Begin with a square $CEFG$. Let a segment $QS$ fall vertically from the top of the square ($GF$) to the bottom of the square ($CE$) and let a segment $CP$ rotate about the point $C$ from the left-hand edge of the square ($CG$) to the bottom of the square ($CE$). Both segments must move at a constant rate and start and end at the same instant. Their intersection traces out the quadratrix, which is the curve joining $G$ and $D$.\footnote{A troublesome detail is that $D$ is not a well-defined point of intersection, but is a limit point of the curve of intersections.} It turns out that $CG/CD=\pi/2$. Thus it is possible to square and rectify a circle using the quadratrix.


%


Despite all of this---Archimedes's theorem relating a circle's circumference to its area, the rectifications of the circle using the spiral and the quadratrix, and the common-sense notion that curves can be straightened like string---Aristotle's belief that curves and line segments can not be compared persisted.
Most famously, in his 1637 \emph{G\'eom\'etrie} Ren\'e Descartes (1596--1650) wrote:
\begin{quote}
Geometry should not include lines [curves] that are like strings, in that they are sometimes straight and sometimes curved, since the ratios between straight and curved lines are not known, and I believe cannot be discovered by human minds, and therefore no conclusion based upon such ratios can be accepted as rigorous and exact. \cite[p.~91]{Descartes:1954}
\end{quote} 

The key to unraveling this mystery is to recognize what it meant for mathematicians to ``know'' the length of a curve. Mathematicians of every era were interested in rectifying curves, but the rules kept changing. It meant that one must construct, using the mathematical tools of the day, a line segment with the same length as the curve (and as algebra overtook geometry it meant being able to express the length of the curve algebraically).  

To the Greeks a curve was rectifiable if it was possible to construct a line segment with the same length as the curve using only compass and straightedge. Thus the  rectifications of the circle using a spiral and a quadratrix did not count as true rectifications. 

They did not count for Descartes either. He expanded the collection of geometric tools beyond lines and circles to include algebraic curves. But he argued that curves like the spiral and the quadratrix were not algebraic. He called such curves ``mechanical,'' because drawing them required the coordination of linear and circular motion. He believed that they were not algebraic because they could be used to square the circle, a feat he believed was impossible  (a proof of this impossibility was still 250 years in the future). In a letter to Marin Mersenne (1588--1648) in 1638 he wrote that 
\begin{quote}
It is against the geometer's style to put forward problems that they cannot solve themselves. Moreover, some problems are impossible, like the quadrature of the circle, etc. \cite{Mancosu:2011}  
\end{quote}

The lack of progress toward curve rectifications gave further evidence that lines and curves were incomparable. In contrast, there were many examples of areas that could be algebraically computed (dating back to Hippocrates's work on lunes \cite[pp.~1--26]{Dunham:1990} and Archimedes's work on the parabola \cite[pp.~233--252]{Archimedes:2002}). As Traub wrote, ``The lack of success concerning the lengths of curves stood in sad contrast to what had been accomplished concerning areas.'' \cite[p.~61]{Traub:1984}

\begin{figure}[ht]
\includegraphics{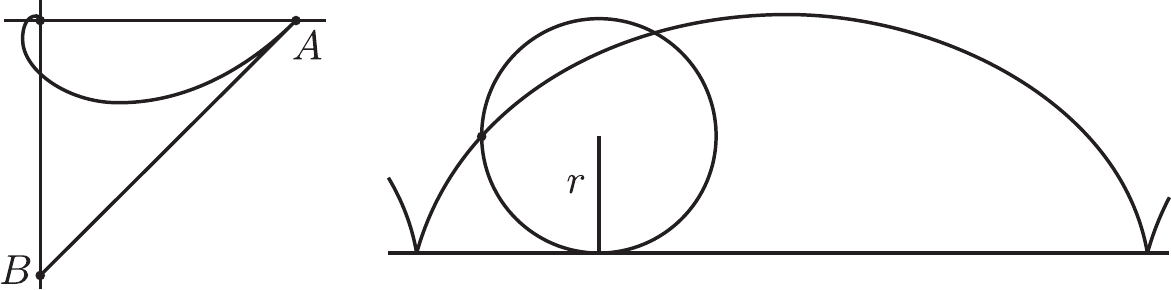}
\caption{The logarithmic spiral and the cycloid}
\label{fig:logspiralcycloid}
\end{figure}

However, there were some steps forward. In 1638 Descartes himself showed that a section of the logarithmic (or equiangular) spiral is Euclidean rectifiable. The portion of  the spiral $r=ae^{\theta}$ for $\theta\le 0$ is the same length as the tangent line segment $AB$  (as in Figure \ref{fig:logspiralcycloid}). This segment is the hypotenuse of an $a\times a$ right triangle, so it is clearly constructible with compass and straightedge. This rectification was rediscovered by Evangelista Torricelli (1608--1647) in 1645, but neither discovery was widely known until 1897 \cite{Boyer:1964}. Recently it was learned that Thomas Harriot (1560--1621) rectified this spiral in the 1590s in the course of his study of map-making and navigation \cite{Pepper:1968}.  

In 1658 Christopher Wren (1632--1723) proved that the cycloid (see Figure \ref{fig:logspiralcycloid}) is Euclidean rectifiable; the length of one arch of the cycloid is equal to eight times the radius of the generating circle. Later Gilles de Roberval (1602--1675) claimed that he had proved this result two decades before Wren. \ncite{Martin:2010}\ncite{Boyer:1964}

But these examples did not strike down Descartes's assertion, for he argued that the logarithmic spiral and the cycloid are not algebraic, but mechanical. They are Euclidean rectifications of transcendental curves. Other mathematicians agreed. In a 1659 letter to Christiaan Huygens (1629--1695), Blaise Pascal (1623--1662) shared this ``beautiful remark'' of Ren\'e de Sluse (1622--1685), 
\begin{quote}
One ought still admire\ldots the order of nature\ldots which does not permit one to find a straight line equal to a curve, except after one has already assumed the equality of a straight line and a curve.'' \cite[p.\ 76]{Traub:1984}
\end{quote}
Sluse and Pascal believed that the argument that the cycloid is rectifiable is circular---it is rectifiable because drawing it requires the coordination of circular motion and straight-line motion.  

\begin{figure}[ht]
\includegraphics{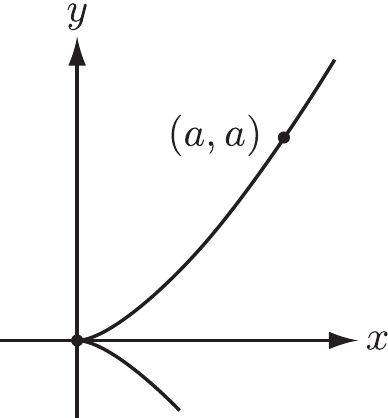}
\caption{}
\label{fig:semicubicalparabola}
\end{figure}

In short, by the end of the 1650s we had examples of algebraic curves that could be rectified using transcendental curves and transcendental curves that could be rectified using algebraic curves. Finally, the impossible happened. In years 1559 and 1660 William Neil (1637--1670) of England, Hendrik van Heuraet (1634--c.\ 1660) of Holland, and Pierre de Fermat (1601--1665) of France each published the same result: an algebraic curve that could be rectified using algebraic (Euclidean, actually) techniques. Fermat wrote of the ``beliefs of the most knowledgable geometers\ldots that it is a law and an order of nature that one can not find a line equal to another curve.'' \cite[p.\ 77]{Traub:1984} Then he presented the semicubical parabola: $ay^{2}=x^{3}$ (see Figure \ref{fig:semicubicalparabola}). The length of this curve from $(0,0)$ to $(a,a)$ is $a(13\sqrt{13}-8)/27$,\ncite{Boyer:1964} and in \emph{G\'eom\'etrie} Descartes proved that if a number can be expressed using the integers, the four arithmetic operations, and square roots, then a segment of that length is constructible with compass and straightedge.

The semicubical parabola finally disproved Descartes's assertion about straight lines and curves. Henk Bos wrote, 
\begin{quote}
The central role of the incomparability of straight and curved in Descartes' geometry was the reason why the first rectifications of algebraic\ldots 
curves in the late 1650s were so revolutionary: they undermined the cornerstone of the edifice of Descartes' geometry. \cite[p.~342]{Bos:2001}
\end{quote}

What about Aristotle's assertion? We have seen several curves that can be rectified using compass and straightedge, so in that sense Aristotle was wrong. However, if we focus not on all curves, but on circles, as Aristotle did in \emph{Physics}, then he was correct. By Archimedes's result we know that to rectify the circle we must be able to construct a segment of length $\pi$ (given a segment of length 1). By Descartes's work we know that a segment is constructible only if its length is an algebraic number, and by Ferdinand von Lindemann's (1852--1939) result from 1882, $\pi$ is transcendental \cite{Lindemann:1882}. Thus it is impossible to rectify the circle. Archimedes's result---relating the area of a circle and its circumference---may have simultaneously proved that $\pi$ is a constant and added more evidence that the circle cannot be rectified, because the problem is equivalent to the quadrature of the circle.

\section{Conclusion}
We would be remiss if we did not mention Archimedes's other contributions to the study of the number $\pi$. We know that $\pi$ has at least two more lives---it is the volume and surface area constants for spheres.  

Just as arc length is trickier than planar area, so too is surface area more subtle than volume. In \emph{Elements} (Prop.~XII.18 \cite[p.~434]{Heath:1908b}) Euclid proved, essentially,  the existence of the volume constant for spheres ($V/d^{3}$). But again it took the genius of Archimedes to prove the that $S/r^{2}$, where $S$ is the surface area of a sphere, is a constant constant independent of the sphere. 

Archimedes did not stop there---he also proved that the circle and sphere constants are intimately related:  \[\pi=\frac{C}{d}=\frac{A}{r^{2}}=6\left( \frac{V}{d^{3}}\right)=\frac{1}{4}\left(\frac{S}{r^{2}}\right).\]  In \emph{On the Sphere and Cylinder} Archimedes brought all of these results together in one concise, elegant result:
\begin{quote}
Every cylinder whose base is the greatest circle in a sphere and whose height is equal to the diameter of the sphere is $3/2$ of the sphere, and its surface together with its bases is $3/2$ of the surface of the sphere. \cite[p.~43]{Archimedes:2002}
\end{quote}
Archimedes recognized this work for the grand achievement that it is. Plutarch (45--120) wrote: 
\begin{quote}
Although [Archimedes] made many excellent discoveries, he is said to have asked his kinsmen and friends to place over the grave where he should be buried a cylinder enclosing a sphere, with an inscription giving the proportion by which the containing solid exceeds the contained. \cite[p.~481]{Plutarch:1917}
\end{quote}

Because of these results about the nature of $\pi$ and his very accurate bounds on the value, it seems fitting that we call $\pi$ \emph{Archimedes's number}.

\bibliographystyle{amsalpha} 
\bibliography{Cd}

\begin{thebibliography}{10}

\bibitem{Archimedes:2002}
Archimedes.
\newblock {\em The works of {A}rchimedes}.
\newblock Dover Publications Inc., Mineola, NY, 2002.
\newblock Reprint of the 1897 edition and the 1912 supplement, Edited by T. L.
  Heath.

\bibitem{Beckmann:1971}
Petr Beckmann.
\newblock {\em A History of Pi}.
\newblock St. Martin's Press, New York, 1971.

\bibitem{Berggren:2004}
Lennart Berggren, Jonathan Borwein, and Peter Borwein.
\newblock {\em Pi, a source book}.
\newblock Springer Verlag, New York, 3 edition, 2004.

\bibitem{Bos:2001}
Henk J.~M. Bos.
\newblock {\em Redefining geometrical exactness}.
\newblock Sources and Studies in the History of Mathematics and Physical
  Sciences. Springer-Verlag, New York, 2001.
\newblock Descartes' transformation of the early modern concept of
  construction.

\bibitem{Boyer:1964}
C.~B. Boyer.
\newblock {\em L'aventure de l'esprit: M\'elanges {A}lexandre {K}oyr\'e},
  volume~I, chapter Early rectification of curves, pages 30--39.
\newblock Hermann, 1964.

\bibitem{Burton:2007}
David~M. Burton.
\newblock {\em The History of Mathematics: An Introduction}.
\newblock McGraw-Hill, Boston, 6 edition, 2007.

\bibitem{Castellanos:1988a}
Dario Castellanos.
\newblock The ubiquitous $\pi$.
\newblock {\em Mathematics Magazine}, 61(2):67--98, 1988.

\bibitem{Castellanos:1988}
Dario Castellanos.
\newblock The ubiquitous $\pi$.
\newblock {\em Mathematics Magazine}, 61(3):148--163, 1988.

\bibitem{Clagett:1964}
Marshall Clagett.
\newblock {\em Archimedes in the {M}iddle {A}ges. {V}ol. {I}: {T}he
  {A}rabo-{L}atin tradition}.
\newblock The University of Wisconsin Press, Madison, Wis., 1964.

\bibitem{Comenetz:2002}
Michael Comenetz.
\newblock {\em Calculus: {T}he elements}.
\newblock World Scientific, River Edge, NJ, 2002.

\bibitem{Descartes:1954}
Ren\'e Descartes.
\newblock {\em The Geometry of {R}en\'e {D}escarte: {T}ranslated from {F}rench
  and {L}atin by {D}avid {E}ugene {S}mith and {M}arcia {L}. {L}atham}.
\newblock Dover Publications Inc., New York, 1954.

\bibitem{Dijksterhuis:1987}
E.~J. Dijksterhuis.
\newblock {\em {A}rchimedes}.
\newblock Princeton University Press, Princeton, NJ, 1987.
\newblock Translated from the Dutch by C. Dikshoorn, Reprint of the 1956
  edition, With a contribution by Wilbur R. Knorr.

\bibitem{Dunham:1990}
William Dunham.
\newblock {\em Journey Through Genius: {T}he Great Theorems of Mathematics}.
\newblock John Wiley \& Sons, New York, 1990.

\bibitem{Gillman:1991}
Leonard Gillman.
\newblock $\pi$ and the limit of $(\sin\alpha)/\alpha$.
\newblock {\em The American Mathematical Monthly}, 98(4):346--349, 1991.

\bibitem{Heath:1980}
Thomas Heath.
\newblock {\em Mathematics in {A}ristotle}.
\newblock Garland Publishing, New York, 1980.

\bibitem{Heath:1908}
Thomas~L. Heath.
\newblock {\em The thirteen books of {E}uclid's {E}lements. {V}olume {I}. Books
  {I}--{II}}.
\newblock The University Press, Cambridge, 1908.

\bibitem{Heath:1908a}
Thomas~L. Heath.
\newblock {\em The thirteen books of {E}uclid's {E}lements. Volume {II}. Books
  {III}--{IX}}, volume~2.
\newblock The University Press, Cambridge, 1908.

\bibitem{Heath:1908b}
Thomas~L. Heath.
\newblock {\em The thirteen books of {E}uclid's {E}lements. Volume {III}. Books
  {X}--{XIII} and appendix}.
\newblock The University Press, Cambridge, 1908.

\bibitem{Jones:1706}
William Jones.
\newblock {\em Synopsis palmariorum mathesios}.
\newblock London, 1706.

\bibitem{Katz:1998}
Victor~J. Katz.
\newblock {\em A history of mathematics: {A}n introduction}.
\newblock Addison Wesley Longman, Reading, MA, 2 edition, 1998.

\bibitem{Knorr:1986}
Wilbur~R. Knorr.
\newblock Archimedes' \emph{{D}imension of the circle}: {A} view of the genesis
  of the extant text.
\newblock {\em Arch. Hist. Exact Sci.}, 35(4):281--324, 1986.

\bibitem{Knorr:1993}
Wilbur~Richard Knorr.
\newblock {\em The ancient tradition of geometric problems}.
\newblock Dover Publications Inc., New York, 1993.
\newblock Corrected reprint of the 1986 original.

\bibitem{Lindemann:1882}
F.~Lindemann.
\newblock \"{U}ber die {Z}ahl $\pi$.
\newblock {\em Mathematische Annalen}, 20(2):213--225, 1882.

\bibitem{Mancosu:2011}
Paolo Mancosu.
\newblock {D}escartes and mathematics.
\newblock In Janet Broughton and John Carriero, editors, {\em A companion to
  {D}escartes}, pages 103--123. Blackwell Publishing, Malden, MA, 2008.

\bibitem{Martin:2010}
John Martin.
\newblock The {H}elen of geometry.
\newblock {\em College Math Journal}, 41(1):17--28, 2010.

\bibitem{Pepper:1968}
Jon~V. Pepper.
\newblock {H}arriot's calculation of the meridional parts as logarithmic
  tangents.
\newblock {\em Arch. History Exact Sci.}, 4(5):359--413, 1968.

\bibitem{Plutarch:1917}
Plutarch.
\newblock {\em {P}lutarch's Lives: {W}ith and {E}nglish translation by
  {B}ernadotte {P}errin}, volume~V.
\newblock William Heinemann, London, 1917.

\bibitem{Richman:1993}
Fred Richman.
\newblock A circular argument.
\newblock {\em The College Mathematics Journal}, 24(2):160--162, 1993.

\bibitem{Ross:1936}
W.~D. Ross.
\newblock {\em {A}ristotle's {P}hysics: {A} revised text with introduction and
  commentary by {W}.\ {D}.\ {R}oss}.
\newblock The Clarendon press, Oxford, 1936.

\bibitem{Schepler:1950a}
Herman~C. Schepler.
\newblock The chronology of pi.
\newblock {\em Mathematics Magazine}, 23(3):165--170, 1950.

\bibitem{Schepler:1950b}
Herman~C. Schepler.
\newblock The chronology of pi.
\newblock {\em Mathematics Magazine}, 23(4):216--228, 1950.

\bibitem{Schepler:1950c}
Herman~C. Schepler.
\newblock The chronology of pi.
\newblock {\em Mathematics Magazine}, 23(5):279--283, 1950.

\bibitem{Seidenberg:1972}
A.~Seidenberg.
\newblock On the area of a semi-circle.
\newblock {\em Arch. History Exact Sci.}, 9(3):171--211, 1972.

\bibitem{Smeur:1970}
A.~J. E.~M. Smeur.
\newblock On the value equivalent to {$\pi$} in ancient mathematical texts. {A}
  new interpretation.
\newblock {\em Arch. History Exact Sci.}, 6(4):249--270, 1970.

\bibitem{Traub:1984}
Gilbert Traub.
\newblock {\em The Development of the Mathematical Analysis of Curve Length
  from {A}rchimedes to {L}ebesgue}.
\newblock PhD thesis, New York University, 1984.

\end{thebibliography}
\end{document}